
\input amssym.def
\input amssym.tex

\def\Titre{Rectangular Scott-type Permanents}

\def\Resume{Soient $x_1, x_2, \ldots, x_n$ les z\'eros d'un polyn\^ome $P(x)$
de degr\'e $n$ et $y_1,y_2, \ldots, y_m$ les z\'eros d'un autre 
polyn\^ome $Q(y)$ de degr\'e $m$. Notre objet d'\'etude est le permanent 
$\per(1/(x_i-y_j))_{1\le i\le n,\, 1\le j\le m}$, 
appel\'e ici permanent de type Scott. Le cas de 
$P(x)=x^n-1$ et $Q(y)=y^n+1$ a \'et\'e consid\'er\'e par R.~F.~Scott. 
Nous pr\'esentons une approche efficace pour d\'eterminer
les \'evaluations explicites des permanents de type Scott, bas\'ee sur des 
g\'en\'eralisations des th\'eor\`emes classiques de Cauchy et Borchardt, 
et d'un th\'eor\`eme r\'ecent de Lascoux. La pr\'esente \'etude prolonge
le travail 
du premier auteur (``G\'en\'eralisation de l'identit\'e de Scott sur les
permanents,'' \`a appara\^\i tre dans Linear Algebra Appl.). 
Notre approche nous permet de fournir de nombreuses
\'evaluations explicites des permanents de type Scott 
pour des choix sp\'eciaux des polyn\^omes $P(x)$ et $Q(y)$, 
y compris des g\'en\'ealisations de tous les r\'esultats de l'article
mentionn\'e ci-dessus et du permanent de Scott lui-m\^eme.
Par exemple, nous prouvons que si $P(x)=x^n-1$ et $Q(y)=y^{2n}+y^n+1$ alors 
le permanent
correspondant de type Scott est \'egal \`a $(-1)^{n+1}n!$. 
} 
\def\Abstract{
Let $x_1,x_2,\dots,x_n$ be the zeroes of a polynomial $P(x)$ of degree $n$ and 
$y_1,y_2,\dots,y_m$ be the zeroes of another polynomial $Q(y)$ of degree $m$.
Our object of study is the permanent $\per(1/(x_i-y_j))_{1\le i\le n,\,
1\le j\le m}$, here named ``Scott-type" permanent, the case of
$P(x)=x^n-1$ and $Q(y)=y^n+1$ having been considered by 
R.~F.~Scott. We present
an efficient approach to determining explicit evaluations of
Scott-type permanents, based on generalizations of
classical theorems by Cauchy and Borchardt, and of a recent
theorem by Lascoux. This continues and extends the work initiated by the first
author (``G\'en\'eralisation de l'identit\'e de Scott sur les
permanents,'' to appear in Linear Algebra Appl.). Our approach enables
us to provide numerous closed form evaluations of Scott-type
permanents for special choices of the polynomials $P(x)$ and
$Q(y)$, including generalizations of all the results from the above
mentioned paper and of Scott's permanent itself. For example, we prove that 
if $P(x)=x^n-1$ and $Q(y)=y^{2n}+y^n+1$ 
then the corresponding Scott-type permanent is equal to $(-1)^{n+1}n!$.
} 


\def\mymatrix#1{\left(\vcenter{\offinterlineskip
     \halign{\vrule height 6pt depth 2pt width 0pt 
      \hfil$\scriptstyle##$\hfil&&\kern
       6pt\hfil$\scriptstyle##$\hfil \crcr#1\crcr}}\right)}

\def\mymatrixserre#1{\left(\vcenter{\offinterlineskip
     \halign{\vrule height 6pt depth 2pt width 0pt 
      \hfil$\scriptstyle##$\hfil&&\kern
       4pt\hfil$\scriptstyle##$\hfil \crcr#1\crcr}}\right)}

\def\per{\hbox{\rm per}}
\def\PER{\hbox{\rm PER}}
\def\D{\Delta}

\def\Res{\hbox{\rm Res}}
\def\Fes{\hbox{\rm Fes}}
\def\FesA{\hbox{$\widetilde{\rm Fes}$}}
\def\mod{\hbox{\rm mod\ }}
\def\NULL{\hbox{\bf 0}}
\def\ID{\hbox{\bf I}}
\def\I{{\cal I}}
\def\s{\sigma}

\def\om{\omega}
\def\gcd{\hbox{\rm gcd}}
\def\sgn{\hbox{\rm sgn}}
\def\frac#1#2{{#1\over#2}}
\def\dfrac#1#2{{{\displaystyle#1}\over{\displaystyle#2}}}
\def\tfrac#1#2{{\textstyle{#1\over#2}}}
\def\({\left(}
\def\){\right)}

\def\qed{\quad\raise -2pt\hbox{\vrule\vbox to 10pt{\hrule width 4pt
   \vfill\hrule}\vrule}}
\def\pointir{\discretionary{.}{}{.\kern.35em---\kern.7em}\nobreak
   \hskip 0em plus .3em minus .4em }
\def\cf{{\it cf}}
\def\rem#1|{\par\medskip{{\it #1}\pointir}}
\def\dem{\par{\it Proof}\pointir}

\def\section#1{\goodbreak\par\vskip .3cm\centerline{\bf #1}
   \par\nobreak\vskip 3pt }

\long\def\th#1|#2\end{\par\medskip
   {\petcap #1\pointir}{\it #2}\par\smallskip}


\def\rAndrews{1}
\def\rBorch{2}
\def\rCallan{3}
\def\rCauchy{4}
\def\rHanInter{5}
\def\rHanScott{6}
\def\rJou{7}
\def\rKit{8}
\def\rKratBI{9}
\def\rKratBN{10}
\def\rLasSquare{11}
\def\rMac{12}
\def\rMincPer{13}
\def\rMincScott{14}
\def\rMuirAB{15}
\def\rSco{16}
\def\rSvr{17}

\def\eDifCtoB{1}
\def\eDefPER{2}
\def\eExpFXY{3}
\def\eExpFXYA{4}
\def\eTwodiag{5}
\def\eDet{6}
\def\eCrazy{7}
\def\eRH{8}
\def\eCirc{9}
\def\eUnit{10}
\def\eCrazyB{11}
\def\espez{12}
\def\eMatrix{13}
\def\ematrix{14}
\def\eCrazyA{15}

\def\tGeneral{1}
\def\tInvolution{2}
\def\tPrincipal{3}
\def\tPrincipalA{4}
\def\tJou{5}
\def\tTwodiag{6}
\def\tCrazy{7}
\def\tCrazyB{8}
\def\tCrazyA{9}
\def\tPerAB{10}
\def\cA{11}
\def\cAa{12}
\def\cAaa{13}
\def\cAb{14}
\def\cAc{15}
\def\cB{16}
\def\cC{17}
\def\cD{18}
\def\cE{19}
\def\cDA{20}
\def\cEA{21}
\def\cF{22}
\def\cG{23}
\def\cGA{24}
\def\cGa{25}
\def\cH{26}
\def\cI{27}
\def\cJ{28}
\def\cK{29}
\def\cKA{30}
\def\cKB{31}
\def\tPerLB{32}
\def\cL{33}
\def\cM{34}
\def\cN{35}
\def\cO{36}
\def\tPerLE{37}
\def\tPerLC{38}
\def\tPerLD{39}
\def\tSumInv{40}
\def\tSumInvThree{41}
\def\tSumInvFour{42}
\def\tSumInvTwo{43}


\magnification=1200
\hsize=11.25cm
\vsize=18cm
\parskip 0pt
\parindent=12pt
\voffset=1cm
\hoffset=1cm

\catcode'32=9

\font\tenpc=cmcsc10
\font\eightpc=cmcsc8
\font\eightrm=cmr8
\font\eighti=cmmi8
\font\eightsy=cmsy8
\font\eightbf=cmbx8
\font\eighttt=cmtt8
\font\eightit=cmti8
\font\eightsl=cmsl8
\font\sixrm=cmr6
\font\sixi=cmmi6
\font\sixsy=cmsy6
\font\sixbf=cmbx6

\skewchar\eighti='177 \skewchar\sixi='177
\skewchar\eightsy='60 \skewchar\sixsy='60

\catcode`\@=11

\def\tenpoint{%
  \textfont0=\tenrm \scriptfont0=\sevenrm \scriptscriptfont0=\fiverm
  \def\rm{\fam\z@\tenrm}%
  \textfont1=\teni \scriptfont1=\seveni \scriptscriptfont1=\fivei
  \def\oldstyle{\fam\@ne\teni}%
  \textfont2=\tensy \scriptfont2=\sevensy \scriptscriptfont2=\fivesy
  \textfont\itfam=\tenit
  \def\it{\fam\itfam\tenit}%
  \textfont\slfam=\tensl
  \def\sl{\fam\slfam\tensl}%
  \textfont\bffam=\tenbf \scriptfont\bffam=\sevenbf
  \scriptscriptfont\bffam=\fivebf
  \def\bf{\fam\bffam\tenbf}%
  \textfont\ttfam=\tentt
  \def\tt{\fam\ttfam\tentt}%
  \abovedisplayskip=12pt plus 3pt minus 9pt
  \abovedisplayshortskip=0pt plus 3pt
  \belowdisplayskip=12pt plus 3pt minus 9pt
  \belowdisplayshortskip=7pt plus 3pt minus 4pt
  \smallskipamount=3pt plus 1pt minus 1pt
  \medskipamount=6pt plus 2pt minus 2pt
  \bigskipamount=12pt plus 4pt minus 4pt
  \normalbaselineskip=12pt
  \setbox\strutbox=\hbox{\vrule height8.5pt depth3.5pt width0pt}%
  \let\bigf@ntpc=\tenrm \let\smallf@ntpc=\sevenrm
  \let\petcap=\tenpc
  \normalbaselines\rm}

\def\eightpoint{%
  \textfont0=\eightrm \scriptfont0=\sixrm \scriptscriptfont0=\fiverm
  \def\rm{\fam\z@\eightrm}%
  \textfont1=\eighti \scriptfont1=\sixi \scriptscriptfont1=\fivei
  \def\oldstyle{\fam\@ne\eighti}%
  \textfont2=\eightsy \scriptfont2=\sixsy \scriptscriptfont2=\fivesy
  \textfont\itfam=\eightit
  \def\it{\fam\itfam\eightit}%
  \textfont\slfam=\eightsl
  \def\sl{\fam\slfam\eightsl}%
  \textfont\bffam=\eightbf \scriptfont\bffam=\sixbf
  \scriptscriptfont\bffam=\fivebf
  \def\bf{\fam\bffam\eightbf}%
  \textfont\ttfam=\eighttt
  \def\tt{\fam\ttfam\eighttt}%
  \abovedisplayskip=9pt plus 2pt minus 6pt
  \abovedisplayshortskip=0pt plus 2pt
  \belowdisplayskip=9pt plus 2pt minus 6pt
  \belowdisplayshortskip=5pt plus 2pt minus 3pt
  \smallskipamount=2pt plus 1pt minus 1pt
  \medskipamount=4pt plus 2pt minus 1pt
  \bigskipamount=9pt plus 3pt minus 3pt
  \normalbaselineskip=9pt
  \setbox\strutbox=\hbox{\vrule height7pt depth2pt width0pt}%
  \let\bigf@ntpc=\eightrm \let\smallf@ntpc=\sixrm
  \let\petcap=\eightpc
  \normalbaselines\rm}

\tenpoint

{
\par
\vskip 30pt plus 24pt minus 3pt\penalty -1000
\vskip 0pt plus -24pt minus 3pt\penalty -1000
\centerline{\bf\Titre}
\vskip 5pt
\penalty 10000 
}
\bigskip
\centerline{Guo-Niu {\petcap Han} et Christian \petcap Krattenthaler} 
{\vskip 6mm\penalty 10000 }
\vbox{\narrower\eightpoint {\petcap Abstract}\pointir\Abstract}
{\vskip 10pt\penalty 10000 }
\vbox{\narrower\eightpoint {\petcap R\'esum\'e}\pointir\Resume}
{\vskip 15pt\penalty 10000 }

\section{1. Introduction}

In 1881, Scott [\rSco] stated, without proof, the following result:

{\smallskip
\parindent=20pt
\narrower\noindent 
{\it Let $x_1, x_2,\ldots, x_n$ be the zeroes of $x^n-1$ and
$y_1,y_2,\ldots, y_n$ be the zeroes of  
$y^n+1$. Let $A$ be the $n\times n$ matrix $(1/(x_i-y_j))_{1\le 
i,j\le n}$. Then
$$
\per(A)=\cases{\displaystyle
(-1)^{n-1\over 2}\ {
n(1\cdot 3\cdot 5\cdots (n-2))^2\over 2^n}, &if $n$ is odd,\cr
0,&if $n$ is even.
}
$$
}}
 
In 1978, in his monograph {\it Permanents}, Minc [\rMincPer, p.~155]
included this result in a list of conjectures on permanents. Since
then, several proofs have been given [\rCallan, \rKit, \rMincScott, \rSvr],
one of which by Minc himself.
All of these proofs are heavily based on the fact that the 
zeroes of the polynomials
$x^n-1$ and $y^n+1$ can be written in simple explicit terms. Thus,
neither of these proofs extends to the more general problem which is
the subject of this paper: 

{\smallskip
\parindent=20pt
\narrower\noindent 
{\it Given a polynomial $P(x)$ of degree $n$ and
a polynomial $Q(y)$ of degree $m$, let $x_1,x_2,\dots,x_n$ be the
zeroes of $P(x)$ and
$y_1,y_2,\dots,\break y_m$ be the zeroes of $Q(y)$. Evaluate the permanent of
the $n\times m$ matrix $(1/(x_i-y_j))_{1\le i\le n,\,
1\le j\le m}$.

\smallskip
}}

\noindent
As usual, the permanent $\per(A)$ of an $n\times m$ matrix $A$ is defined
as the sum of all possible products of $n$ coefficients of $A$ chosen such that
no two of the coefficients are taken from the same row nor from the
same column (see [\rMincPer, Ch.~1, (1.1)]). Given the assumptions of
the above problem, we call the permanent of the matrix 
$(1/(x_i-y_j))_{1\le i\le n,\,1\le j\le m}$ a {\it Scott-type
permanent}, and denote it by $\PER(P(x),Q(y))$.

\medskip
In [\rHanScott], the first author presented a new approach to this
type of problem in the case $n=m$, i.e., in the case that both
polynomials have the same degree. This approach does not rely at all on
explicit analytic forms of zeroes of polynomials. Instead, it makes essential
use of recent symmetric functions techniques, in particular of a
theorem due to Lascoux [\rLasSquare], which the latter author 
established in his
\'etude on the square ice model of statistical mechanics.

In the present paper, we are going to extend this approach to arbitrary
$n$ and $m$. This requires extensions of classical theorems of Cauchy
and Borchardt (see Theorems (Cauchy+) and (Borchardt+) in Section~2),
and an extension of Lascoux's theorem (see Theorem (Lascoux+)). 
As a result (see Theorem~\tGeneral), we are able to express any
Scott-type permanent as the quotient of a determinant which features
complete homogeneous and elementary symmetric functions in the zeroes
of the two polynomials,
divided by the resultant of the two polynomials. 
In particular, it follows immediately that any Scott-type permanent
is rational in the coefficients of the polynomials $P(x)$ and $Q(y)$.

In Section~5 we apply this result to obtain explicit evaluations of
Scott-type permanents in numerous special cases. Amongst others,
we provide generalizations of all the results from [\rHanScott],
thus also covering Scott's permanent itself. For the proofs of the results in
Section~5, we make use of two particular specializations of our main
theorem, Theorem~\tGeneral, which we derive in Section~3 (see
Theorems~\tPrincipal\ and \tPrincipalA), and of
four determinant evaluations, which we state and establish
separately in Section~4.

Finally, we also comment briefly on an alternative approach to the
evaluation of Scott-type determinants, due to the first author
[\rHanInter]. It allows to express Scott-type permanents in terms of
weighted sums over involutions (see Theorem~\tInvolution\ in Section~2). By
combining this result with some of the evaluations in Section~5, we
obtain interesting summation theorems, which are presented in
Section~6. 


\section{2. The general theory}

In [\rHanScott], the main ingredients are theorems by Cauchy,
Borchardt, Lascoux, and a lemma on the resultant. Since we intend to
extend the approach of [\rHanScott] to the case of {\it
rectangular} Scott-type permanents (corresponding to polynomials of,
possibly, different degrees), we have to first provide the appropriate
extensions of these theorems.

Given positive integers $m$ and $n$ and two sets 
$X=\{x_1,x_2,\dots,
x_n\}$ and $Y=\{y_1,y_2,\dots,y_m\}$ of variables, we
use the following notations:
$$R(X,Y):=\prod_{i=1}^n \prod_{j=1}^m (x_i-y_j) 
\hbox{\quad and \quad }
\D(X):=\prod_{i<j}(x_i-x_j).$$
$$
{
\def\CauchyUn#1#2{\displaystyle { 1\over {x_#1}-y_#2} }
\def\Cauchy#1{\CauchyUn{#1}1 & \CauchyUn{#1}2 &\cdots &\CauchyUn{#1}m}
 \Bigl({1\over x_i-y_j}\Bigr):=
\pmatrix{
\Cauchy{1}\cr
\Cauchy{2}\cr
\vdots&\vdots&\vdots&\vdots\cr
\noalign{\vskip 3pt}
\Cauchy{n}\cr
}.
}
$$

We first state a variation on Cauchy's evaluation of his double alternant
(\cf. [\rCauchy; \rMuirAB, vol.~1, pp.~342--345]).

\th Theorem (Cauchy+)|
For $m\geq n$, let
$$
{
\def\CauchyUn#1#2{\displaystyle { 1\over {x_#1}-y_#2} }
\def\Cauchy#1{\CauchyUn{#1}1 & \CauchyUn{#1}2 &\cdots &\CauchyUn{#1}m}
C(X,Y):=
\pmatrix{
\Cauchy{1}\cr
\Cauchy{2}\cr
\vdots&\vdots&\vdots&\vdots\cr
\noalign{\vskip 3pt}
\Cauchy{n}\cr
\noalign{\vskip 3pt}
1 & 1 & \cdots & 1\cr
y_1 & y_2 & \cdots & y_m\cr
y_1^2 & y_2^2 & \cdots & y_m^2\cr
\vdots&\vdots&\vdots&\vdots\cr
y_1^{m-n-1} & y_2^{m-n-1} & \cdots & y_m^{m-n-1}\cr
}.
}
$$
Then
{
\def\CauchyUn#1#2{\displaystyle { 1\over {x_#1}-y_#2} }
\def\Cauchy#1{\CauchyUn{#1}1 & \CauchyUn{#1}2 &\cdots &\CauchyUn{#1}m}
$$
\det(C(X,Y))
=(-1)^{n(n-1)/2} {\D(X)\D(Y)\over R(X,Y)}.
$$
}
\end
\dem
If $n=m$, this is exactly Cauchy's theorem. The general case can be
either established directly, or, it may be observed that
the ``general" case is in fact {\it implied\/} by Cauchy's theorem. To see
this, consider the above identity with $n=m$. Given $k<m$, expand both sides as
power series in $1/x_{k+1}$, \dots, $1/x_m$, and compare coefficients
of $1/x_{k+1}x_{k+2}^2\cdots x_m^{m-k}$ on both sides.
\qed

Next we state the required extension of Borchardt's theorem [\rBorch;
\rMuirAB, vol.~2, pp.~173--175]. It can be established
by reading through the proof of Borchardt's theorem given
in [\rAndrews, Proof of Cor.~5.1], 
ignoring however the restriction $m=n$ (see
also [\rHanInter]). 

\th Theorem (Borchardt+)|
For $m\geq n$, let
$$
{
\def\CauchyUn#1#2{\displaystyle { 1\over ({x_#1}-y_#2)^2} }
\def\Cauchy#1{\CauchyUn{#1}1 & \CauchyUn{#1}2 &\cdots &\CauchyUn{#1}m}
B(X,Y):=
\pmatrix{
\Cauchy{1}\cr
\Cauchy{2}\cr
\vdots&\vdots&\vdots&\vdots\cr
\noalign{\vskip 3pt}
\Cauchy{n}\cr
\noalign{\vskip 3pt}
1 & 1 & \cdots & 1\cr
y_1 & y_2 & \cdots & y_m\cr
y_1^2 & y_2^2 & \cdots & y_m^2\cr
\vdots&\vdots&\vdots&\vdots\cr
y_1^{m-n-1} & y_2^{m-n-1} & \cdots & y_m^{m-n-1}\cr
}.
}
$$
Then
$$
\det(B(X,Y))=\det(C(X,Y)) \times \per\({1\over x_i-y_j}\).
$$
\end

Our next goal is to derive the required extension of (a special case
of) Lascoux's theorem
[\rLasSquare, Theorem~q].
Let $Z=\{z_1,z_2,\dots,z_n\}$ be another set of variables, of equal
cardinality as $X$.
For $1\leq i\leq n$ we define the divided difference 
$\partial_i$ by
$$
\partial_i: f\mapsto {f-f^{\s_i}\over x_i-z_i},
$$
where $\s_i$ is the transposition which interchanges $x_i$ and $z_i$.
It is easy to see that
$$
\partial_i {1\over z_i-y} =  {1\over(x_i-y)(z_i-y)}.
$$
Since the operator
$\partial_i$ acts only on one row  
of the matrix $C(Z,Y)$ (to be precise, the $i$-th row), it follows that
$$
\partial_1\partial_2\cdots\partial_n
(\det(C(Z,Y))|_{Z=X}=\det(B(X,Y)).\eqno{(\eDifCtoB)}
$$
Now, to generalize Lascoux's theorem to our case, one simply
reads\break 
through
the proof of Theorem~q in [\rLasSquare], on introducing slight
modifications if necessary. The result is:

\th Theorem (Lascoux+)|
Let
$H(X)$ be the $n\times(m+n-1)$ matrix defined by
$$H(X):=\big(h_{j-i}(X)\big)_{1\leq i\leq n, 1\leq j\leq
m+n-1},$$
where $h_s(X)$ denotes the complete homogeneous symmetric function of
degree $s$ in the variables $X$ (\cf. [\rMac, Ch.~1]),
and $E(Y)$ be the $(m+n-1)\times n$ matrix defined by
$$
E(Y):=\Bigl(
(j-2k+2)(-1)^{m-j+k-1}e_{m-j+k-1}(Y)
\Bigr)_{1\leq j\leq m+n-1, 1\leq k \leq n},
$$
where $e_s(Y)$ denotes the elementary symmetric function of
degree $s$ in the variables $Y$ (\cf. [\rMac, Ch.~1]).
Then
$$
\partial_1\partial_2\cdots\partial_n
(\D(Z)R(X,Y))|_{Z=X}=\D(X) \det\bigl(H(X)E(Y)\bigr).
$$
\end

Now we are in the position to state our main theorem, which will
enable us, in Section~5, to evaluate numerous Scott-type permanents
in closed form. The theorem implies immediately that any Scott-type
permanent $\PER(P(x),Q(y))$ 
is a rational function in the coefficients of the two
polynomials $P(x)$ and $Q(y)$.

\th Theorem \tGeneral|
Let $m$ and $n$ be arbitrary positive integers, and let
$X=\{x_1,x_2,\dots,x_n\}$ and $Y=\{y_1,y_2,\dots,y_m\}$ be two sets of
variables. Then
$$\per\({1\over x_i-y_j}\)={\det\bigl(H(X)E(Y)\bigr)\over R(X,Y)}.$$
\end

\dem
First let $m\ge n$.
By combining Theorems (Cauchy+) and (Borchardt+), and Equation
(\eDifCtoB), we obtain
$$(-1)^{n(n-1)/2}\frac {R(X,Y)} {\Delta(X)\Delta(Y)}
\partial_1\partial_2\cdots\partial_n
(\det(C(Z,Y))|_{Z=X}$$
for the permanent.
Now we apply Theorem (Cauchy+) again in order to replace $C(Z,Y)$ by
the corresponding product form guaranteed by the theorem. After having
used that the divided differences $\partial_i$ commute with
$\Delta(Y)/R(X,Y)R(Z,Y)$ (because the latter expression is symmetric
in $x_i$ and $z_i$), Theorem (Lascoux+) applies and yields the desired
result.

If $m<n$, the permanent clearly vanishes.
According to Theorem (Lascoux+), it suffices to establish that
$$
U:=\partial_1\partial_2\cdots\partial_n(\D(Z)R(X,Y))=0.
$$
To begin with, we rewrite the Vandermonde determinant evaluation as
$\D(Z)R(X,Y)=(-1)^{n(n-1)/2}\det(z_i^{j-1}R(x_i,Y))$. Writing
$Q(x_i):=\break R(x_i, Y)=a_mx_i^m+\cdots+a_1x_i+a_0$,
we obtain for $U$ the expression
$$
U=\det(\partial_i z_i^{j-1}R(x_i,Y))=
\det\bigl( z_i^{j-1}Q(x_i) -
x_i^{j-1}Q(z_i)\bigr)\Big/\prod(x_i-z_i).
$$
Since for $m<n$, we have
$$
\sum_{j=0}^{m} a_m \bigl(z_i^{j}Q(x_i) - x_i^{j}Q(z_i)\bigr)=0,
$$
the $m+1$ elements $ z_i^{j}Q(x_i) - x_i^{j}Q(z_i)$, $0\leq j\leq m$,
are linearly dependent. Hence, $U=0$.
\qed

\medskip
In [\rHanInter], the first author obtained another expression for
the permanent, in form of a certain weighted sum over involutions.
To state and explain this formula, for $s\in X$ define
$$
L(s;X,Y):=\sum_{x\not=s}{1\over x-s} + \sum_{y\in Y}{1\over s-y}.
$$

Let us denote by $\I(n)$ the set of involutions on $\{1,2,\dots,n\}$.
Given an involution $\s\in\I(n)$, we define the weight $\Psi(\s)$ of
$\s$ by
$$\Psi(\s;X,Y):= 
\prod_{(ij)\in \s}{1\over
(x_i-x_j)^2} \prod_{(k)\in\s}
{ L(x_k;X,Y)}
$$
where the first product is over all transpositions $(ij)$ in the
disjoint cycle decomposition of $\s$, and where the second product is
over all fixed points $k$ of $\s$. Then the result from [\rHanInter]
is the following.
\medskip

\th Theorem \tInvolution|
Let $m$ and $n$ be arbitrary positive integers, and let $X=\{x_1,x_2,\dots, 
x_n\}$ and $Y=\{y_1,y_2,\dots,y_m\}$ be two sets of variables. Then
$$
\per\({1\over x_i-y_j}\)=\sum_{\s\in\I(n)}\Psi(\s;X,Y).
$$
\end

\rem Example 1|
Let $n=1$, $X=\{x\}$,  $Y=\{y_1,y_2,\dots,y_m\}$. Then
$$
\eqalign{
H(X)&:=\big(h_{j-1}(X)\big)_{i=1, 1\leq j\leq m}
=(h_0(x), h_1(x), \ldots, h_{m-1}(x)),\cr
E(Y)&:=\Bigl(
j(-1)^{m-j}e_{m-j}(Y)
\Bigr)_{1\leq j\leq m, k=1}\cr
&=(1(-1)^{m-1} e_{m-1}(Y), 2(-1)^{m-2} e_{m-2}(Y),\dots, m e_0(Y))^t.
\cr
}
$$
We have
$$
H(X)E(Y)=\sum_{j=0}^{m-1} (j+1)(-1)^{m-j-1}e_{m-j-1}(Y)x^j.
$$
Therefore,
$$
\eqalign{
\per\({1\over x-y}\)
&={1\over x-y_1} +{1\over x-y_2} +\cdots +{1\over x-y_m}  
\quad\hbox{\rm [Definition, Th.~\tInvolution]}\cr 
&={\sum_{j=0}^{m-1}(j+1)(-1)^{m-j-1}e_{m-j-1}(Y)x^j\over 
(x-y_1)(x-y_2)\cdots(x-y_m)}.
\quad\hbox{\rm [Th.~\tGeneral]}\cr 
}
$$

\medskip

\rem Example 2|
For $n=4$ and $m=3$, Theorem~\tGeneral\ yields the following identity:
$$
\let\thematrice=\mymatrix
\det\(\thematrice{
h_0&h_1&h_2&h_3&h_4&h_5\cr
0    &h_0&h_1&h_2&h_3&h_4\cr
0  &0  &  h_0&h_1&h_2&h_3\cr
0  &0  &0  &   h_0&h_1&h_2\cr
}\times
\thematrice{
e_2 & e_3 & 0 & 0 \cr
-2e_1 & 0 & 2e_3 & 0\cr
3e_0 & -e_1 & -e_2 & 3e_3 \cr
0 & 2e_0 & 0 & -2e_2 \cr
0 & 0 & e_0 & e_1 \cr
0 & 0 & 0 & 0 \cr
}\)=0.
$$
For $n=2$ and $m=1$, Theorem~\tInvolution\ yields
$$
{1\over (x_1-x_2)^2} + 
\({1\over x_2-x_1}+{1\over x_1-y}\)
\({1\over x_1-x_2}+{1\over x_2-y}\)=0.
$$

\section{3. The case of $P(x)=x^n-1$ and of $P(x)=x^{n-1}+\cdots+x+1$}

In this section we specialize Theorem~\tGeneral\ to the case that
the $x_i$'s are the zeroes of the
polynomial $P(x)=x^n-1$ or of $P(x)+x^{n-1}+\cdots+x+1$, 
and the $y_i$'s are the zeroes of an
arbitrary other polynomial. (This covers, for example, the
case of Scott's identity). For the remainder of this section,
we fix $m$ and $n$, $m\ge n$.

Let $X=\{x_1, x_2, \ldots, x_n\}$ be the set of zeroes of $x^n-1$, and
let $Y=\{y_1, y_2, \ldots,
y_m\}$ be the set of zeroes of
$Q(x)=a_my^m+a_{m-1}y^{m-1}+\cdots+a_1y^1+a_0$, with $a_m=1$. 
We write
$$
\PER(P,Q):={\per\({1\over
x_i-y_j}\)}={\det(H(X)E(Y))\over R(X,Y)}.
\eqno{(\eDefPER)}
$$
Since 
$$\sum h_i(X) t^i = {1\over \prod_i(1-tx_i)}
= {1\over (1-t^{n})}=1+t^n+t^{2n}+\cdots,$$ 
we have
$$
h_k(X)=\cases{
1, &if $k=0\ (\mod n)$,\cr
0, &if $k\not=0\ (\mod n)$.\cr
}
$$
We denote by $\ID_k$ the $k\times k$ identity matrix, and by
$\NULL_{l,c}$ the $l\times c$ matrix with all entries equal to $0$.
For all $r$, we write $r\%n$ for the number between $1$ and $n$ that
satisfies 
$r\ (\mod n)=r\%n\ (\mod n)$. Then we have
$$
{
\def\oneIn{\ID_n\ \vrule}
H(X)=\pmatrix{
\oneIn&\oneIn&\cdots\ \vrule&\oneIn&
\matrix{\ID_{m'}\cr \NULL_{n-{m'},{m'}}\cr
}\cr
},
}
$$
with $m'=(m+n-1)\%n$.

Furthermore, let ${\rm diag}_n^i(c_1, c_2, \ldots, c_n)$ denote the
$n\times n$ ``diagonal" matrix, in which the (broken) diagonal starts
in the $i$-th row,
$$
\pmatrix{
\NULL_{i-1, n-i+1} & {\rm diag}_{i-1}(c_{n-i+2},c_{n-i+3},\dots, 
c_{n})
\cr {\rm diag}_{n-i+1}(c_1, c_2,\dots, c_{n-i+1}) & 
\NULL_{n-i+1, i-1} \cr
}.
$$
According to the definition of $E(Y)$, a simple calculation yields that
$$
H(X)E(Y)=
\sum_{r=0}^m {\rm diag}_n^{r\%n}
(ra_r, (r-1)a_r, \dots, (r-n+1)a_r).
$$

\rem Example 3|
For $n=3$ and $m=4$, let $P(x)=x^3-1$ and $Q(y)=y^4+a_3y^3+a_2y^2+a_1y+a_0$.
Then $H(X)E(Y)$ is the sum of the following matrices:
$$
\let\thematrice=\mymatrix
\thematrice{
0&-a_0&0\cr
0&0&-2a_0\cr
0&0&0\cr
}+
\thematrice{
a_1&0&0\cr
0&0&0\cr
0&0&-a_1\cr
}+
\thematrice{
0&0&0\cr
2a_2&0&0\cr
0&a_2&0\cr
}+
\thematrice{
0&2a_3&0\cr
0&0&a_3\cr
3a_3&0&0\cr
}+
\thematrice{
4&0&0\cr
0&3&0\cr
0&0&2\cr
}.
$$
\medskip

In the theorem below, we summarize our findings.

\th Theorem \tPrincipal|
Let $P(x)=x^n-1$ and
$Q(y)=a_my^m+\cdots+a_1y^1+a_0$, $a_m$ not necessarily $1$.
Writing
$$
\Fes(Q)= \det
\biggl(\sum_{r=0}^m {\rm diag}_n^{r\%n}
(ra_r, (r-1)a_r, \dots, (r-n+1)a_r)
\biggr),\eqno{(\eExpFXY)}
$$
we have
$$
\PER(P, Q)={\Fes(Q)\over \Res(P,Q)},
$$
where $\Res$ is the classical resultant of two polynomials.
\end
\dem
We have already seen that the theorem is true for $a_m=1$. On the
other hand, there hold
$
\Fes(\lambda Q)=\lambda^n \Fes(Q)$
and
$\Res(P, \lambda Q)=\lambda^n \Res(P,Q)
$, as is easily verified.
\qed

\medskip

Now let us consider the case that
$X=\{x_1,x_2,\dots,x_{n-1}\}$ is the set of zeroes of
$P(x)=x^{n-1}+\cdots+x+1$. We perform an analysis very similar to the
one before, using the fact that we have
$$\sum h_i(X) t^i = {1\over \prod_i(1-tx_i)}
= {1-t\over \prod_i(1-t^n)}=1-t+t^n-t^{n+1}+-\cdots.$$ 
In order to state the result, we introduce the following notation: We
write $\widetilde{\rm diag}{}_{n-1}^i(c_1, c_2, \ldots, c_{n-1})$ for the
$(n-1)\times(n-1)$ matrix
$$
\pmatrix{
\NULL_{i-1, n-i} &\NULL_{i-1,1}& {\rm
diag}_{i-2}(c_{n-i+2},c_{n-i+3},\dots,  c_{n-1}) 
\cr 
{\rm diag}_{n-i}(c_1, c_2,\dots, c_{n-i}) &\NULL_{n-i,1}& 
\NULL_{n-i+1, i-2} \cr
}
$$
if $i>1$, respectively ${\rm diag}_{n-1}(c_1, c_2,\dots, c_{n-1})$ if $i=1$.
This is again a matrix with a (possibly broken) diagonal, in which the
diagonal ``jumps over" one row and column in the case that it is
broken. (Note
the slight discrepancy in dimension between the diagonal and the zero
matrices.)

\th Theorem \tPrincipalA|
Let $P(x)=x^{n-1}+\cdots+x+1$ and
$Q(y)=a_my^m+\cdots+a_1y^1+a_0$, $a_m$ not necessarily $1$.
Writing
$$
\eqalign{
\FesA(Q)=& \det
\biggl(\sum_{r=0}^m \widetilde{\rm diag}{}_{n-1}^{r\%n}
(ra_r, (r-1)a_r, \dots, (r-n+2)a_r)
\cr
&- \sum_{r=0}^m \widetilde{\rm diag}{}_{n-1}^{(r-1)\%n}
(ra_r, (r-1)a_r, \dots, (r-n+2)a_r)
\biggr),
}\eqno{(\eExpFXYA)}
$$
we have
$$
\PER(P, Q)={\FesA(Q)\over \Res(P,Q)},
$$
where, again, $\Res$ is the classical resultant of two polynomials.
\end

\medskip

According to Theorems~\tPrincipal\ and \tPrincipalA, 
for accomplishing the evaluation of the permanent,
it is necessary to evaluate the numerator
$\Fes(Q)$, respectively $\FesA(Q)$, and the denominator $\Res(P,Q)$.
For the evaluation of the resultant $\Res(P,Q)$, we make use of the following
lemma, which, for example, appears explicitly in [\rJou]. (In fact, it
follows from a special case of Proposition~\tTwodiag\ in the next section.)

\th Lemma \tJou|Let $d$ be the greatest common divisor of $m$ and
$n$, and let $A$ and $C$ be two nonzero constants.
Then the resultant of the two polynomials
$Ax^m-B$ and $Cx^n-D$ is given by
$$
\Res(Ax^m-B, Cx^n-D)=(-1)^m
\bigl(A^{n/d}D^{m/d}-B^{n/d}C^{m/d}
\bigr)^d.
$$
\end

What concerns the evaluation of $\Fes(Q)$, respectively $\FesA(Q)$, 
we refer the reader to the
next section for the determinant evaluations that we are going to
use. In combination with Lemma~\tJou, these
will allow us to evaluate Scott-type permanents in numerous
special cases, see Section~5.

\section{4. Determinant evaluations}

\th Proposition \tTwodiag|
Let $n$ and $r$ be positive integers, $r\le n$, and
$x_1,x_2,\dots,x_n$, $y_1,y_2,\dots,y_n$ be indeterminates. Then, with
$d=\gcd(r,n)$, we have
$$
\eqalign{&
\det\pmatrix{
x_1&0&\dots&0&y_{n-r+1}&0&\cr
0&x_2&0&&0&y_{n-r+2}&0\cr
&&\ddots&&&&\ddots&0\cr
0&&&&&&0&y_n\cr
y_1&0&\cr
0&y_2&0\cr
&0&\ddots&0&&&\ddots&0\cr
&&0&y_{n-r}&0&&0&x_n\cr}\cr
&\hskip3.6cm=\prod _{i=1} ^{d}\bigg(\prod _{j=1} ^{n/d}x_{i+(j-1)d}-
(-1)^{n/d}\prod _{j=1} ^{n/d}y_{i+(j-1)d}\bigg).
}\eqno{(\eTwodiag)}
$$
(I.e., in the matrix there are only nonzero entries along two
diagonals, one of which is a broken diagonal.)
\end

\dem
Let $A=(A_{ij})_{1\le i,j\le n}$ be the matrix of which the
determinant is taken in (\eTwodiag). By definition, we have 
$$\det
A=\sum _{\s\in S_n} ^{}\sgn\, \s \prod _{i=1} ^{n}A_{i\s(i)},
\eqno(\eDet)
$$ 
where
$S_n$ denotes the symmetric group of order $n$. Since the matrix $A$
is a sparse matrix, only those
permutations $\s$ contribute to the sum which have the property
$\s(i)=i$ or $\s(i)=i+r$ mod $n$ for all $i$. Thus, the decomposition
into disjoint cycles of such a permutation consists only of cycles
of length $1$, and of cycles of length $n/d$ of the form
$(i,i+r,\dots,i+(n/d-1)r)$ (where, again, all integers have to be taken
modulo $n$). Using these observations and the fact that the sign of
any cycle of length $n/d$ is $(-1)^{n/d-1}$ 
in (\eDet) yields (\eTwodiag) immediately. \qed

\th Theorem \tCrazy|
Let $a,b,c,d,e$ be indeterminates. For any positive integer $n$ and
integers $i,j$ let $n(i,j)$ denote 1 plus
(the representative between $0$ and $n-1$ of)
the residue class of $i-j$ mod $n$.
Then
$$
\displaylines{
\det_{1\le i,j\le n}\big((n(i,j)+c)(n(i,j)a+b)+d
-(j-1)(n(i,j)a+e)\big)\kern1cm\cr
\eightpoint
=\det
\pmatrix{
(c+1)(a+b)\hfill&(c+n)(na+b)\hfill&&(c+2)(2a+b)\hfill\cr
+d&+d-(na+e)&\dots&+d-(n-1)(2a+e)\cr
\noalign{\vskip 3pt}
(c+2)(2a+b)&(c+1)(a+b)\hfill&&(c+3)(3a+b)\hfill\cr
+d&+d-(a+e)&\dots&+d-(n-1)(3a+e)\cr
\vdots&\vdots&\dots&\vdots\cr
\noalign{\vskip 3pt}
(c+n)(na+b)&(c+n-1)\hfill&&(c+1)(a+b)\hfill\cr
+d&\times((n-1)a+b)+d&\dots&+d-(n-1)(a+e)\cr
&\hfill-((n-1)a+e)&&\cr}\cr
\hfill\hfill=(-n)^{n-1}U_n(a,b,c,d,e)\prod _{i=3} ^{n}(ia+b+ca),
\hfill(\eCrazy)}
$$ 
where $U_n(a,b,c,d,e)$ is the polynomial
$$
\eightpoint
\eqalign{
&U_n(a,b,c,d,e)=\frac {(n+1)(n+2)} {3} a^2 + 
  \frac {(n+1)(2n+7)} {6} a b + 
  \frac {(n+1)} {2} b^2\cr
&\quad  +
   \frac {(n+1)(2n+7)} {6}a^2c + 
   \frac {(3n+5)} {2} a b c + 
  b^2c + \frac {(n+1)} {2} a^2c^2 + a b c^2
   + \frac {(n+3)} {2} a d \cr
&\quad + b d + a c d - 
  \frac {(n-1)(2n+5)} {6} a e - \frac {(n-1)} {2} b e -
   \frac {(n-1)} {2} a c e.
}
$$
In the case that $n=2$, the product in {\rm(\eCrazy)} has to be read as $1$,
and in the case that $n=1$, the product has to be interpreted as $1/(2a+b+ca)$.
\end
\dem
In the cases $n=1$ and $n=2$, the claim can be verified directly. 
For $n\ge3$, 
we use the ``identification of factors'' method as explained in
[\rKratBN, Sec.~2.4] or [\rKratBI, Sec.~2].

We proceed in several steps. An outline is as
follows. In the first step we show that $\prod _{i=3} ^{n}(ia+b+ca)$ 
is a factor of the determinant
as a polynomial in $a,b,c,d,e$. In the second step we prove that
$U_n(a,b,c,d,e)$ is a factor of the determinant.
Then, in the third step, 
we determine the maximal degree of the determinant as a
polynomial in $a$, and also in $b$, $c$, $d$, and in $e$.
It turns out that the maximal degree is $n$ as a polynomial in $a$,
the same being true as a polynomial in $b$ and as a polynomial in $c$,
while it is 1 as a polynomial in $d$, the same being true as a
polynomial in $e$.
On the other hand, the degree in $a$, and also in $b$ and in $c$, 
of the product on
the right-hand side of (\eCrazy), which by the first two steps divides
the determinant, is exactly $n$. It is exactly 1 in $d$ and also in
$e$. Therefore we are forced to conclude that
the determinant equals
$$C(n)U_n(a,b,c,d,e)\prod _{i=3} ^{n}(ia+b+ca),\eqno(\eRH)$$
where $C(n)$ is a constant independent of $a,b,c,d,e$.
Finally, in the fourth step, we determine the constant $C(n)$, which
turns out to equal $(-n)^{n-1}$.
Clearly, this would finish the proof of theorem.

\smallskip
{\it Step 1. For $i=3,\dots,n$ the term $(ia+b+ca)$ is a factor of
the determinant}. We claim that, if $b=-ia-ca$, we have
$$
\eqalign{
&-\big(\hbox{row }(n-i)\big)+3\big(\hbox{row }(n-i+1)\big)\cr
&\hskip3cm
-3\big(\hbox{row }(n-i+2)\big)+\big(\hbox{row }(n-i+3)\big)=0
}
$$
as long as $n>i$. (Here, (row $i$) denotes the $i$-th row of the
matrix underlying the determinant in (\eCrazy).) 
In the case that $n=i$, we claim that we have
$$
-3\big(\hbox{row }1\big)+3\big(\hbox{row }2\big)
-\big(\hbox{row }3\big)+\big(\hbox{row }n\big)=0
$$
as long as $n>3$, and that we have
$
-\big(\hbox{row }1\big)+\big(\hbox{row }2\big)=0
$
if $n=3$. All these claims are easily verified by an obvious
case-by-case analysis.

\smallskip
{\it Step 2. The polynomial $U_n(a,b,c,d,e)$ is a factor of
the determinant}. We claim that if $d$ is chosen so that
$U_n(a,b,c,d,e)$ vanishes, we have   
$$\sum _{j=1} ^{n}((j+1)a+b+ca)\big(\hbox{column }j\big)=0.$$
Again, it is a routine task to verify this identity.

\smallskip
{\it Step 3. The determinant is a polynomial in $a$ (in $b$,
respectively in $c$) of maximal degree
$n$, and a polynomial in $d$ (respectively in $e$) of maximal degree
$1$}.  The first claim follows from the fact that
each term in the defining expansion
of the determinant has degree $n$ in $a$ (as well as in $b$,
respectively in $c$). To establish the second claim, we simply
subtract the first row of the determinant 
from all other rows, with the effect that only the entries in the first
row contain $d$ and $e$ after these transformations.
Since the right-hand side of (\eCrazy), which
by Steps~1 and 2 divides the determinant as a polynomial in $a,b,c,d,e$,
also has degree $n$ in $a$, in $b$, and in $c$, and degree $1$ in $d$, and
in $e$, the determinant 
and the right-hand side of
(\eCrazy) differ only by a multiplicative constant.

\smallskip
{\it Step 4. 
The evaluation of the multiplicative constant.} By the preceding
steps we know that the determinant equals (\eRH). 
In particular, if we set $a=c=d=e=0$
and $b=1$, we have
$$\det_{1\le i,j\le n}\(n(i,j)\)
=C(n)(n+1)/ {2}.\eqno(\eCirc)$$
The matrix on the left-hand side of (\eCirc) is a circulant
matrix with entries $1,2,\dots,n$. Hence, its determinant equals
$$\prod _{\om\hbox{\sevenrm \ : zero of }x^n-1}
^{}(1+2\om+3\om^2+\cdots+n\om^{n-1}).$$
The sum is easily evaluated by observing that it is
the derivative of a geometric series. It turns out to be equal to
$-n/(1-\om)$. The
resulting product simplifies by the observation 
$$\prod _{\om\hbox{\sevenrm \ : zero of }x^n-1,\ \om\ne1} ^{}(1-\om)=
(1+x+\cdots+x^{n-1})\big|_{x=1}=n.\eqno(\eUnit)$$
Thus, the
determinant in (\eCirc) equals $(-n)^{n-1}{(n+1)/ 2}$. Therefore
$C(n)$ is equal to $(-n)^{n-1}$.

\medskip
This finishes the proof of (\eCrazy) and thus of the theorem.
\qed

\th Theorem \tCrazyB|
Let $a$ be an indeterminate. For any positive integer $n$ and
integers $i,j$ let $s(i,j)$ denote 
(the representative between $0$ and $n-1$ of)
the residue class of $i-j+1$ mod $n$.
Then
$$
\displaylines{
\det_{1\le i,j\le n-1}\kern-3pt\(\eightpoint
\left\{\matrix {(n-m-1) + j(1+a-n) 
&
\hbox{if }i=j-2\ (\mod n)\cr
\noalign{\vskip 3pt}
(n-1)(m-1) + j(1-a) 
&
\hbox{if }i=j-3\ (\mod n)\cr
\noalign{\vskip 3pt}
(n-m-3-2s(i,j)) +j
&\hbox{otherwise}\hfill}\right\}\)\cr
\hfill\hfill=
(-1)^{n-1}\frac {1} {n}\prod _{i=2} ^{n}(nm-ia).
\hfill(\eCrazyB)}
$$ 
\end
\dem
The matrix underlying this determinant is a matrix whose elements have
a uniform definition, except for two (broken) diagonals, the one with 
$i=j-2\ (\mod n)$, and the one with $i=j-3\ (\mod n)$.

To begin with, we reorder the rows so that the next-to-last row
becomes the first row, the last row becomes the second row, and then
follow the remaining rows in their original order.

Now we add all the rows to the (new) last row. In the resulting
matrix, we change the sign of the last row and, subsequently, move it 
up so that it becomes the third
row. As a result of these manipulations, we obtain
$$
(-1)^{n-1}\kern-7pt\det_{1\le i,j\le n-1}
\(\left\{
\matrix{
(n-1) (m-1) + j (1-a)&
\hbox{if }i=j\cr
(n-m-1) + j (a-n+1)&
\hbox{if }i=j+1\cr
-m-n+1-2 i+3 j&
\hbox{if }i<j\cr
-m+n+1-2 i+3j&
\hbox{if }i>j+1\cr
}
\right\}\).
\eqno{(\espez)}
$$

Next we apply further row operations. We subtract the second row from
the first, the third from the second, \dots, the $(n-1)$-st row from the
$(n-2)$-nd row, in that order. Subsequently, we repeat the same kind
of operations, but stop before the last row, i.e., we subtract the
second row from 
the first, the third from the second, \dots, the $(n-2)$-nd row from the
$(n-3)$-rd row, in that order. As a result, the above determinant is
converted into the determinant of the following matrix
$$
\(\left\{
\matrix{
nm-ja
&\hbox{if }i=j-2\hfill\cr
-2nm+2n-j(3a-n)
&\hbox{if }i=j-1\hbox{ and }i< n-2\hfill\cr
nm-2n-j(3a-2n)
&\hbox{if }i=j\hbox{ and }i<n-2\hfill\cr
j(a-n)
&\hbox{if }i=j+1\hbox{ and }i<n-2\hfill\cr
2
&\hbox{if }i=n-2\hbox{ and }j\le n-4\hfill\cr
(n-3)(a-n)+2
&\hbox{if }i=n-2\hbox{ and }j= n-3\hfill\cr
nm-2n+2 -(n-2)(2a-n)
&\hbox{if }i=n-2\hbox{ and }j= n-2\hfill\cr
n-nm+1+(n-1)(a-1)
&\hbox{if }i=n-2\hbox{ and }j= n-1\hfill\cr
-m-n+3+3j
&\hbox{if }i=n-1\hbox{ and }j\le n-3\hfill\cr
n-m-1+(n-2)(a-n+1)
&\hbox{if }i=n-1\hbox{ and }j= n-2\hfill\cr
(n-1)(m-a)
&\hbox{if }i=n-1\hbox{ and }j= n-1\hfill\cr
0
&\hbox{otherwise}\hfill\cr
}
\right\}\).
\eqno{(\eMatrix)}
$$
This is a matrix with four ``special'' diagonals (the diagonals with
$i=j-2$, $j-1$, $j$, $j+1$, respectively) and two ``special" rows (the last two
rows). All other entries are zero.

Now we factor $(nm-3a)$ out of the first row. Subsequently, we
substract $(a-n)$ times the (new) first row from the second. Now one
is able to factor $(nm-4a)$ out of the (new) second row. Next, we 
substract $2(a-n)$ times the (new) second row from the third row. Etc.
We stop this procedure in the $(n-3)$-rd row. Thus, the determinant of
the matrix in (\eMatrix) is equal to $\prod _{i=3} ^{n-1}(nm-ia)$ times the
determinant of the matrix
$$
\pmatrix{
1&-2&\hphantom{-}1&0&\dots\cr
0&\hphantom{-}1&-2&1&0&\dots\cr
\vdots&\ddots&\ddots&\ddots&\ddots&&\vdots\cr
\vdots&&\ddots&\ddots&\ddots&\ddots&0\cr
&&&0&1&-2&1\cr
a_{n-2,1}&a_{n-2,2}&\multispan4{\quad \hbox to3cm{\dotfill}}&a_{n-2,n-1}\cr
a_{n-1,1}&a_{n-1,2}&\multispan4{\quad \hbox to3cm{\dotfill}}&a_{n-1,n-1}\cr
},
\eqno{(\ematrix)}
$$
where the entries in the last two rows are still the same as in
(\eMatrix). 

We now perform the final set of transformations. We add column~1
through column $n-2$ to column $n-1$, and then we add
$$\sum _{j=1} ^{n-3}(n-j-1)\big(\hbox{column }j\big)$$
to column $n-2$. The effect is that a block matrix is obtained of the
form 
$$\pmatrix{U&0\cr *&M},$$
where $U$ is an $(n-3)\times (n-3)$ upper-triangular matrix with 1s on
the diagonal, and where $M$ is the $2\times 2$ matrix
$$\pmatrix{
nm-n-2a+2&n-2\cr
\tfrac {1} {2}(n-2)(n+m+2a-nm-3)&\tfrac {5} {2}n-\tfrac {n^2} {2}
+m-a-3
}.
$$
Clearly, the determinant of $U$ is 1, while the determinant of $M$ is
$(m-a)(nm-2a)$. Putting everything together, we have completed the
proof of (\eCrazyB).
\qed

\th Corollary \tCrazyA|
Let $a$ be an indeterminate. For any positive integer $n$ and
integers $i,j$ let $s(i,j)$ denote 
(the representative between $0$ and $n-1$ of)
the residue class of $i-j+1$ mod $n$.
Then
$$
\displaylines{
\det_{1\le i,j\le n-1}\(\left\{\matrix {-2s(i,j)-a
+j-1&\hbox{if }i\ne j-2\ (\mod n)\cr
(n-1)(n+a-j-1)&\hbox{if }i=j-2\ (\mod n)}\right\}\)\cr
\hfill\hfill=n^{n-2}\prod _{i=0} ^{n-2}(i+a).
\hfill(\eCrazyA)}
$$ 
\end
\dem
In the determinant, we move the last row on top, replace the (now)
last row by the sum of all the rows, factor $(-1)$ out of the
resulting row, and finally move it
up so that it becomes the second row, retaining the order of all the
other rows. These operations did not change
the value of the determinant. However, the resulting determinant is
exactly the determinant in (\espez) with $a=n$ and $m=a+n$. 
As we have shown in the proof of Theorem~\tCrazyB, the latter
determinant differs from the determinant in (\eCrazyB) just by a sign
of $(-1)^{n-1}$. Thus, we
obtain the right-hand side of (\eCrazyA).
\qed

\section{5. Closed form evaluations for Scott-type permanents}

\th Theorem \tPerAB|
Let $n$, $m$ and $r$ be positive integers  and $d=\gcd(n,r)$. Then
$$
\eqalign{
&\PER\bigg(x^n-1, \sum _{\ell =0} ^{m}a_\ell y^{\ell n}+\sum _{\ell =0}
^{m}b_\ell y^{\ell n+r}\bigg)\cr
&\quad =
-\dfrac {
\eqalign{
&d^n\prod _{i=1} ^{d}\Bigg(\bigg(\sum _{\ell =0} ^{m}a_\ell
\bigg)^{n/d}\bigg(
\sum _{\ell =0} ^{m}(i-n\ell-1) a_\ell 
\bigg/\sum _{\ell =0} ^{m}da_\ell  \bigg)_{n/d}\cr
&\hskip1cm-\bigg(-\sum _{\ell =0} ^{m}b_\ell \bigg)^{n/d}\bigg(
\sum _{\ell =0} ^{m}(i-r-n\ell-1) b_\ell 
\bigg/\sum _{\ell =0} ^{m}db_\ell  \bigg)_{n/d}\Bigg)
}
} {\Bigg(\bigg(\sum _{\ell =0} ^{m}a_\ell \bigg)^{n/d}-
\bigg(-\sum _{\ell =0} ^{m}b_\ell \bigg)^{n/d}\Bigg)^d},
}
$$
where $(\alpha)_k$ is the standard notation for {\it shifted factorials},
$(\alpha)_k:=\alpha(\alpha+1)\cdots(\alpha+k-1)$, $k\ge1$, and $(\alpha)_0:=1$.
\end
\dem
Let us first consider the case that $n\nmid r$.
According to Theorem~\tPrincipal, we have to compute the quotient
$\Fes(Q)/\Res(x^n-1,Q)$, where  $Q=\sum _{\ell =0} ^{m}a_\ell y^{\ell n}+\sum
_{\ell =0} ^{m}b_\ell y^{\ell n+r}$. In order to compute $\Fes(Q)$,  
in (\eExpFXY) we replace $a_{\ell n}$ by $a_\ell $ and $a_{\ell n+r}$ by
$b_\ell $,
$\ell =0,1,\dots,m$, and set all other $a_\ell $'s equal to zero. In the
resulting determinant we move the last row to the top, thus creating a
sign of $(-1)^{n-1}$, and finally apply
Proposition~\tTwodiag\ with $x_j=\sum _{\ell =0} ^{m}(n\ell-j+1) a_\ell$
and $y_j=\sum _{\ell =0}
^{m}(n\ell+r-j+1) b_\ell $, $j=1,2,\dots,n$.
For the evaluation of
$\Res(x^n-1,Q)$ we note that
$$
\eqalign{
&\Res\bigg(x^n-1,\sum _{\ell =0} ^{m}a_\ell y^{\ell n}+\sum _{\ell =0}
^{m}b_\ell y^{\ell n+r}\bigg)\cr
& =\kern -3pt\prod _{\om\hbox{\sevenrm \ : zero of }x^n-1} ^{}\bigg(
\sum _{\ell =0} ^{m}a_\ell +\om^r\sum _{\ell =0}
^{m}b_\ell \bigg)
=\Res\bigg(x^n-1,y^r\sum _{\ell =0}
^{m}b_\ell+\sum _{\ell =0} ^{m}a_\ell\bigg).
}
$$
Now we can apply Lemma~\tJou.

If on the other hand $n\mid r$, then $d=n$. It can be verified
directly that the claimed formula remains valid in that case, too.
\qed

\rem Remark|
Given a polynomial $Q(y)$, there is no unique way to write it in the form
$ \sum _{\ell =0} ^{m}a_\ell y^{\ell n}+\sum _{\ell =0}
^{m}b_\ell y^{\ell n+r}$. For example, we may write
$Q(y)=y^n+a+b$ as $Q(y)=(y^n+ay^0)+b$, or as $Q(y)=(y^n+(a+b))$ (i.e.,
either with $a_1=1$, $a_0=a$, $r=0$, $b_0=b$, or with 
$a_1=1$, $a_0=a+b$, $b_\ell=0$ for all $\ell$).
Regardless which choice we make, Theorem~\tPerAB\ yields
$$
\PER(x^n-1, y^n+a+b)=(-1)^{n+1}
{\prod_{i=1}^{n} (i-(n-i)(a+b))\over (a+b+1)^n}.
$$

\th Corollary \cA|
Let $n$ and $m$ be positive integers. Then
$$
\PER\bigg(x^n-1, \sum _{\ell =0} ^{m}a_\ell y^{\ell n}\bigg)=
- {\bigg(-n\sum _{\ell =0} ^{m}\ell a_\ell \bigg/
\sum _{\ell =0} ^{m}a_\ell  \bigg)_n}.
$$
\end

\th Corollary \cAa|
We have
$$
\PER(x^n-1, y^{mn}+\cdots+y^{2n}+y^{n}+1)={
-(-mn/2)_n
}.
$$
\end

\th Corollary \cAaa|
If $m$ is even then
$$
\PER(x^n+1, y^{mn}+\cdots+y^{2n}+y^{n}+1)={
(-mn/2)_n
}.
$$
\end
\dem
We use the case $a_\ell=(-1)^\ell$ in Corollary~\cA, and the fact that
$$\per\(\frac {1} {x_i\root n\of{-1}-y_j}\)_{1\le i\le n,\, 1\le j\le
mn}
\kern-5pt=
-\per\(\frac {1} {x_i-y_j/\root n\of{-1}}\)_{1\le i\le n,\, 1\le j\le mn}.$$
\qed

\th Corollary \cAb|
We have
$$
\PER\bigg(x^n-1, \sum _{\ell=0} ^{m}\ell y^{\ell n}\bigg)={
-(-n(2m+1)/3)_n
}.
$$
\end

\th Corollary \cAc|
We have
$$
\PER\bigg(x^n-1, \sum _{\ell=0} ^{m}\ell y^{\ell^2 n}\bigg)={
-(-nm(m+1)/2)_n
}.
$$
\end

\th Corollary \cB|
We have
$$
\PER(x^n-1, y^{mn}+ay^{rn}+b)={
-(-(m+ra)n/(a+b+1))_n
}.
$$
\end

\th Corollary \cC|
We have
$$
\PER(x^n-1, y^{mn}+1)={
-(-mn/2)_n
}.
$$
\end

For $m=1$ one recovers Scott's identity stated at the beginning of the
introduction.

\th Corollary \cD|
If $m+ra=a+b+1\ne0$, then
$$
\PER(x^n-1, y^{mn}+ay^{rn}+b)={
(-1)^{n+1} n!
}.
$$
\end

\th Corollary \cE|
If $a\ne-2$ then
$$
\PER(x^n-1, y^{2n}+ay^{n}+1)={
(-1)^{n+1} n!
}.
$$
\end

\th Corollary \cDA|
If $m+ra=0$ and $a+b+1\ne0$, then
$$
\PER(x^n-1, y^{mn}+ay^{rn}+b)={
0
}.
$$
\end

\th Corollary \cEA|
If $b\ne1$ then
$$
\PER(x^n-1, y^{2n}-2y^{n}+b)={
0
}.
$$
\end

\th Corollary \cF|
Let $n$ and $m$ be positive integers and $d=\break \gcd(n,m)$. Then
$$
\PER(x^n-1, y^m+b)=
-\dfrac {d^n\prod _{i=1} ^{d}\(\(\tfrac {i-m-1} {d}\)_{n/d}-
(-b)^{n/d}\(\tfrac {i-1} {d}\)_{n/d}\)} {\(1-(-b)^{n/d}\)^d}.
$$
\end
\dem
In Theorem~\tPerAB, set $m=0$, $a_0=b$, $b_0=1$, and replace $r$ by
$m$, in this order.
\qed

\th Corollary \cG|
If $\gcd(m,n)=1$, then
$$
\PER(x^n-1, y^m+b)=(-1)^{n+1} {
m(m-1)\cdots(m-n+1) \over 1-(-b)^n
}.
$$
\end

\th Corollary \cGA|
If $\gcd(m,n)=1$, then
$$
\eqalign{
&\PER(x^{s(n-1)}+\cdots+x^{2s}+x^s+1, y^{s(m-1)}+\cdots+y^{2s}+y^s+1)\cr
&\hskip3cm
= {
\frac {\displaystyle
\prod _{i=0} ^{s-1}\bigg(\prod _{\ell=0} ^{n-1}(i+\ell s)-\prod
_{\ell=0} ^{n-1}(i+\ell s-ms)\bigg)} {(mns)^s}
}.
}
$$
\end
\dem
Consider the Scott-type permanent $\PER(x^{sn}-1, y^{sm}-q^{sm})$ 
(expressed in
terms of its definition). When we multiply it by $(1-q)^s$ and then
perform the limit $q\to 1$, then the permanent reduces to
$\PER(x^{s(n-1)}+\cdots+x^s+1, y^{s(m-1)}+\cdots+y^s+1)$, 
as is straightforward to see. On the other hand, the permanent that we
started with is the permanent in Corollary~\cF\ with $n$ replaced by
$sn$, $m$ replaced by $sm$, and $b=-q^{sm}$. Indeed,
if we multiply the right-hand side from Corollary~\cF\ (with these
choices for the parameters)
by $(1-q)^{sm}$, and then perform the limit $q\to1$, we obtain
exactly the claimed result.
\qed

\th Corollary \cGa|
If $\gcd(m,n)=1$, then
$$
\eqalign{
&\PER(x^{n-1}+\cdots+x+1, y^{m-1}+\cdots+y+1)\cr
&\hskip5cm
=(-1)^{n+1} {
(m-1)\cdots(m-n+1) \over n
}.
}
$$
\end

\th Corollary \cH|
If $\gcd(m,n)=1$ and $n$ is odd, then
$$
\PER(x^n-1, y^m+1)= {
m(m-1)\cdots(m-n+1) \over 2
}.
$$
\end

\th Corollary \cI|
If $n$ is odd, then
$$
\PER(x^n-1, y^{n+1}+1)= {
(n+1)! \over 2
}.
$$
\end

\th Corollary \cJ|
Let $n$ and $r$ be positive integers (not necessarily $n>r$) 
and $d=\gcd(n,r)$. Then
$${
\def\nd{{n\over d}}
\eqalign{
&\PER(x^n-1, y^n+ay^r+b)\cr
&\quad =
-\dfrac {
d^n\prod _{i=1}^{d}\biggl(
(b+1)^\nd
\Bigl({ib-b+i-n-1\over d(b+1)}\Bigr)_\nd
-(-a)^\nd
\Bigl({i-r-1\over d}\Bigr)_\nd
\biggr)}
{\Bigl((b+1 )^\nd-(-a)^\nd\Bigr)^d}.
}}
$$
\end
\dem
In Theorem~\tPerAB, set $m=1$, $a_0=b$, $a_1=1$, $b_0=a$, and $b_1=0$.
\qed

\th Corollary \cK|
Let $n$ and $r$ be positive integers (not necessarily $n>r$) 
and $\gcd(n,r)=1$. Then
$${
\PER(x^n-1, y^n+ay^r+b) =
(-1)^{n+1}\ \frac {
\prod_{i=1}^{n}(i-(n-i)b))
-a^n
(-r)_n
}
{(b+1 )^n-(-a)^n}.
}
$$
\end

If $1\leq r\leq n-1$ then $(-r)_n=0$. Thus, one recovers the results in
[\rHanScott]. On the other hand, if we set $r=n+1$, we obtain, for
example, the following result.

\th Corollary \cKA|
We have
$${
\PER(x^n-1, y^{n+1}+y^n-1) =
n^n-(-1)^n(n+1)!.
}
$$
\end

\th Corollary \cKB|
We have
$${
\PER(x^n-1, y^{n}+ny-1) =
1.
}
$$
\end

\th Theorem \tPerLB|
Let $n$ and $m$ be positive integers and $a$ be an arbitrary
number. Then
$$
\eqalign{
&\PER\bigg(x^n-1, \sum _{\ell =0} ^{mn-1}(\ell +a)y^{\ell }\bigg)\cr
&\hskip2cm=
(-1)^{n-1}\frac {n(m-1)\,V_n(a,m)} {6\({mn} +2a-1\)}\( {a}
+(m-1)n+1\)_{n-2}  ,
}
$$
where $V_n(a,m)$ is the polynomial
$$
V_n(a,m)=1 - 6  a + 6 a^2 +  n - 2 a n -
      5  m n
    +  10 a m n -  m n^2 + 4  m^2 n^2.
$$
\end
\dem
Let $Q=\sum _{\ell =0} ^{mn-1}(\ell +a)y^{\ell }$.
Using
Theorem~\tPrincipal\ again, we have to compute 
$\Fes(Q)/\Res(x^n-1,Q)$. 
In order to compute $\Fes(Q)$, in (\eExpFXY) set $a_{\ell }=\ell+a$,
$\ell =0,1,\dots,mn-1$. In the
resulting determinant we move the last row to the top, thus creating a
sign of $(-1)^{n-1}$, and finally apply
Theorem~\tCrazy\ with $a=1$, $b=(m-1)n-1$, $c=a-1$, 
$d=n^2(m-1)(2m-1)/6-an(m-1)/2$ and
$e=a+n(m-1)/2-1$. 

For the computation of $\Res(x^n-1,Q)$ we note that
$$
\eqalign{
&\Res(x^n-1, Q)=
\prod _{\om\hbox{\sevenrm \ : zero of }x^n-1} ^{}\bigg(
\sum _{\ell =0} ^{mn-1}(\ell +a)\om^\ell \bigg)\cr
&\hskip1.8cm
=\({mn\choose 2}+mna\)\prod _{\om\hbox{\sevenrm \ : zero of }x^n-1,\
\om\ne1} ^{}\bigg(
\sum _{\ell =1} ^{mn}\ell\om^{\ell-1} \bigg).
}
$$
The sum in the last line is the derivative of 
a geometric series, and is therefore
easily evaluated. The result of the summation 
turns out to be 
$-mn/(1-\om)$. The computation is completed by the observation (\eUnit),
and some simplification.
\qed
\medskip

\th Corollary \cL|
Let $n$ and $m$ be positive integers, $n\ge2$. Then
$$
\PER\bigg(x^n-1, \sum _{\ell =0} ^{mn-1}\ell y^{\ell }\bigg)
=(-1)^{n-1} {(4mn-n-1) (mn-2)! \over 6 (mn-n-1)!}.
$$
\end

\th Corollary \cM|
Let $n$ and $m$ be positive integers, $n\ge2$. Then
$$
\PER\bigg(x^n-1, \kern-2pt\sum _{\ell =0} ^{mn-1}(\ell+1) y^{\ell }\bigg)
=(-1)^{n-1} {(4mn-n+1) (mn-n) (mn-1)! \over 6 (mn-n+1)!}.
$$
\end

\th Corollary \cN|
Let $n$ and $m$ be positive integers, $n\ge2$. Then
$$
\PER\bigg(x^n-1, \sum _{\ell =0} ^{mn-1}(mn-\ell) y^{\ell }\bigg)
={(m-1)\,(n+1)! \over 6 }.
$$
\end

\th Corollary \cO|
Let $n$ and $m$ be positive integers, $n\ge2$. Then
$$
\PER\bigg(x^n-1, \sum _{\ell =0} ^{mn-1}(mn-\ell-1) y^{\ell }\bigg)
= {(m-1)\,n! \over 6 }.
$$
\end

\rem Remark|
It is also possible to move forward and derive formulas for
$\PER\(x^n-1, \sum _{\ell =0} ^{mn/s-1}(\ell +a)y^{\ell s}\)$,
where $s$ is some positive integer. This would require to find
analogues of Theorem~\tCrazy\ in which $n(i,j)$ is replaced by
$D\, n(i,j)$, where $D$ is the inverse of $s/\gcd(n,s)$ modulo
$n/\gcd(n,s)$. As calculations aided by the
computer indicate, the resulting determinant evaluations have forms
very similar to (\eCrazy). That is, the result shows a product of
linear factors in $a$ and $b$ as the one on the right-hand side of
(\eCrazy), and one irreducible polynomial of higher degree (such as
$U_n(a,b,c,d,e)$ in (\eCrazy)). However, as $D$ increases, the
degree of the irreducible polynomial also increases, whereas the
amount of linear factors in $a$ and $b$ decreases, so that the results
become more and more unwieldy. We therefore content ourselves with
stating the result when $s$ divides $n$.

\th Theorem \tPerLE|
Let $n$, $m$ and $s$ be positive integers so that $s\mid n$, 
and let $a$ be an arbitrary
number. Then
$$
\eqalign{
&\PER\bigg(x^n-1, \sum _{\ell =0} ^{\frac{mn}s-1}(\ell +a)y^{\ell s}\bigg)
=
(-1)^{n-1}\frac {s^{n-2s}} {6^s\,(m n+2a s-s)^s}\cr
&\hskip2cm
\times  
\prod _{k=0} ^{s-1}\(\(a+\tfrac {1} {s}(n m-n-k)+1\)_{n/s-2}V_{n,s}(a,m,k)\)
,
}
$$
where $V_{n,s}(a,m,k)$ is the polynomial
$$
\eightpoint
\eqalign{
V_{n,s}&(a,m,k)=
6 k^2 m n + 6 k m n^2 - 10 k m^2 n^2 + m n^3 - 5 m^2 n^3 + 4 m^3 n^3
 -  6 k^2 s\cr
& + 12 a k^2 s - 6 k n s + 12 a k n s + 12 k m n s - 
    24 a k m n s - n^2 s
 + 2 a n^2 s + 6 m n^2 s \cr
&- 12 a m n^2 s - 
    5 m^2 n^2 s + 10 a m^2 n^2 s - 2 k s^2 + 12 a k s^2 - 12 a^2 k s^2 - 
    n s^2\cr
& + 6 a n s^2 - 6 a^2 n s^2 + m n s^2 - 6 a m n s^2 + 6 a^2 m n s^2.
}
$$
\end
\dem
Let $Q=\sum _{\ell =0} ^{mn-1}(\ell +a)y^{\ell s}$.
Again, according to
Theorem~\tPrincipal, we have to compute 
$\Fes(Q)/\Res(x^n-1,Q)$. 
In order to compute $\Fes(Q)$, in (\eExpFXY) set $a_{\ell s}=\ell+a$,
$\ell =0,1,\dots,mn/s-1$, and all other $a_i$'s to zero. In the
resulting determinant we move the last row to the top, thus creating a
sign of $(-1)^{n-1}$. Since only every $s$-th $a_i$ is nonzero, we are
dealing with a determinant of a matrix 
in which a lot of entries are zero. If we
permute rows and columns so that first come the rows and columns whose
indices are congruent $1$ mod $s$, then come the rows and columns whose
indices are congruent $2$ mod $s$, etc., then the matrix of which we
want to compute the determinant assumes a block form, with 
$(n/s)\times (n/s)$ blocks on the diagonal, and zeroes
otherwise. Therefore the determinant equals the product of the
determinants of the $s$ matrices of dimension $(n/s)\times (n/s)$ on the
diagonal. Each of these determinants can be evaluated by means of 
Theorem~\tCrazy. We leave it to the reader to fill in the details.

For the computation of $\Res(x^n-1,Q)$ we proceed as in the proof of
Theorem~\tPerLB. 
\qed

\th Theorem \tPerLC|
Let $n$ and $m$ be positive integers, $n\ge2$, and $a$ be an arbitrary
number. Then
$$
\eqalign{
&\PER\bigg(x^{n-1}+\cdots+x+1, \sum _{\ell =0} ^{mn-1}(\ell +a)y^{\ell }\bigg)\cr
&\hskip4cm=
(-1)^{n-1}
\( {a}
+(m-1)n+1\)_{n-1}  .
}
$$
\end
\dem
Let $P(x)=x^{n-1}+\cdots+x+1$
and $Q(y)=\sum _{\ell =0} ^{mn-1}(\ell +a)y^{\ell }$.
This time we apply
Theorem~\tPrincipalA. According to that theorem, we have to compute 
$\FesA(Q)/\Res(P,Q)$. 
In order to compute $\FesA(Q)$, in (\eExpFXYA) we set $a_{\ell }=\ell+a$,
$\ell =0,1,\dots,mn-1$. The resulting determinant is 
exactly $m^{n-1}$ times the determinant in
Corollary~\tCrazyA\ with $a$ 
replaced by $((m-1)n+a+1)$. For the
computation of the resultant of
$P$ and $Q$ we proceed as in the proof of
Theorem~\tCrazy. Simplification of the result yields the claimed
expression. 
\qed

\th Theorem \tPerLD|
Let $n$ and $m$ be positive integers, $n\ge2$, and $a$ be an arbitrary
number. Then
$$
\eqalign{
&\PER\bigg(x^{n-1}+\cdots+x+1, \sum _{\ell =0} ^{mn-2}(\ell +a)y^{\ell }\bigg)\cr
&\hskip4cm=
(-1)^{n-1}\frac {(nm-n)_{n-1}\,(nm+a-1)^{n-1}} 
{(mn+a-1)^n-(a-1)^n}.
}
$$
\end
\dem
Let $P(x)=x^{n-1}+\cdots+x+1$
and $Q(y)=\sum _{\ell =0} ^{mn-1}(\ell +a)y^{\ell }$.
Using Theorem~\tPrincipalA\ again, we have to compute 
$\FesA(Q)/\Res(P,Q)$. 
In order to compute $\FesA(Q)$, in (\eExpFXYA) we set $a_{\ell }=\ell+a$,
$\ell =0,1,\dots,mn-2$. The resulting determinant is 
exactly $m^{n-1}$ times the determinant in
Theorem~\tCrazyB\ with $m$ replaced by $mn+A-1$, $a$ replaced by
$(mn+A-1)/m$, and $A$ replaced
by $a$, in that order.

For the computation of $\Res(x^n+\cdots+x+1,Q)$ 
we proceed similarly as in the proof of Theorem~\tPerLB. Using an
observation from that proof, we note that
$$
\eqalign{
&\Res(x^n+\cdots+x+1, Q)=
\prod _{\om\hbox{\sevenrm \ : zero of }x^n-1,\ x\ne1} ^{}\bigg(
\sum _{\ell =0} ^{mn-2}(\ell +a)\om^\ell \bigg)\cr
&\hskip1.2cm
=\prod _{\om\hbox{\sevenrm \ : zero of }x^n-1,\ x\ne1} ^{}
\(-\frac {mn} {1-\om}-(mn+a-1)\om^{mn-1}\)\cr
&\hskip1.2cm
=(-1)^{n-1}\prod _{\om\hbox{\sevenrm \ : zero of }x^n-1,\ x\ne1} ^{}
\frac {mn+a-1+\om(1-a)} {\om(1-\om)}\cr
&\hskip1.2cm
=\frac {(mn+a-1)^n-(a-1)^n} {n^2m}.
}
$$
The result follows now upon some simplification.
\qed

\section{6. Sums of involutions}

As in [\rHanScott, Sec.~4], we may obtain interesting summation theorems by
combining special evaluations of Scott-type determinants (see
Section~5) with Theorem~\tInvolution. For example, if we combine
Corollary~\cE\ with Theorem~\tInvolution, we obtain the following
result. 

\th Proposition \tSumInv|Let
$x_1,x_2,\ldots, x_n$ be the zeroes of
$x^n-1$. Then
$$
\sum_{\s\in\I(n)} \  \prod_{(ij)\in \s}{1\over
(x_i-x_j)^2} \prod_{(k)\in\s}
{n+1\over 2 x_k}
=(-1)^{n+1}n!,
$$
where the first product is over all transpositions $(ij)$ of $\s$ in its
disjoint cycle decomposition, and where the second product is over all 
fixed points $k$ of $\s$.
\end

Similarly, if we combine Corollary~\cEA\ with Theorem~\tInvolution, we
obtain the result below.
\th Proposition \tSumInvThree|Let 
$x_1,x_2,\ldots, x_n$ be the zeroes of
$x^n-1$. Then
$$
\sum_{\s\in\I(n)} \  \prod_{(ij)\in \s}{1\over
(x_i-x_j)^2} \prod_{(k)\in\s}
{n-1
\over 
2 x_k}
={0},
$$
where the first product is over all transpositions $(ij)$ of $\s$ in its
disjoint cycle decomposition, and where the second product is over all 
fixed points $k$ of $\s$.
\end

If we combine Corollary~\cKB\ with Theorem~\tInvolution, then we
obtain the following result.
\th Proposition \tSumInvFour|Let 
$x_1,x_2,\ldots, x_n$ be the zeroes of
$x^n-1$. Then
$$
\sum_{\s\in\I(n)} \  \prod_{(ij)\in \s}{1\over
(x_i-x_j)^2} \prod_{(k)\in\s}
{2n+(n+1)x_k
\over 
2 x_k^2}
={1},
$$
where the first product is over all transpositions $(ij)$ of $\s$ in its
disjoint cycle decomposition, and where the second product is over all 
fixed points $k$ of $\s$.
\end

Finally, if we combine Corollary~\cI\ with Theorem~\tInvolution, we
obtain the result below.
\th Proposition \tSumInvTwo|Let $n$ be odd, and let
$x_1,x_2,\ldots, x_n$ be the zeroes of
$x^n-1$. Then
$$
\sum_{\s\in\I(n)} \  \prod_{(ij)\in \s}{1\over
(x_i-x_j)^2} \prod_{(k)\in\s}
{1-n+(3+n)x_k
\over 
2 (1+x_k)x_k}
={(n+1)!\over 2},
$$
where the first product is over all transpositions $(ij)$ of $\s$ in its
disjoint cycle decomposition, and where the second product is over all 
fixed points $k$ of $\s$.
\end

\def\article#1|#2|#3|#4|#5|#6|#7|
    {{\leftskip=7mm\noindent
     \hangindent=2mm\hangafter=1
     \llap{[#1]\hskip.35em}{\petcap#2}\pointir
     #3, {\sl #4} {\bf #5} ({\oldstyle #6}),
     pp.\nobreak\ #7.\par}}
\def\livre#1|#2|#3|#4|
    {{\leftskip=7mm\noindent
    \hangindent=2mm\hangafter=1
    \llap{[#1]\hskip.35em}{\petcap#2}\pointir
    {\sl #3}\pointir #4.\par}}
\def\divers#1|#2|#3|
    {{\leftskip=7mm\noindent
    \hangindent=2mm\hangafter=1
     \llap{[#1]\hskip.35em}{\petcap#2}\pointir
     #3.\par}}

\vfill\break\vglue 1cm
\vskip 1cm
\centerline{\bf References}\bigskip\bigskip
{
\eightpoint

\article \rAndrews|G. E. Andrews, I. P. Goulden, D. M. Jackson|Generalizations of
Cauchy's summation theorem for Schur functions|Trans. Amer. Math.
Soc.|310|1988|805--820|

\article \rBorch|C. W. Borchardt|Bestimmung der symmetrischen 
Verbindungen vermittelst ihrer erzeugenden Funktion|Crelle J.|53|1855|193--198|

\article \rCallan|D. Callan|On evaluating permanents and a matrix of 
cotangents|Linear and Multilinear Algebra|38|1995|193--205|

\article \rCauchy|A. L. Cauchy|M\'emoire sur les fonctions altern\'ees
et les sommes altern\'ees|Exercices d'analyse 
et de phys. math.|2|1841|151--159|

\divers \rHanInter|G.-N. Han|Interpolation entre Cauchy et Borchardt, in
preparation, {\oldstyle 1999}|

\divers \rHanScott|G.-N. Han|G\'en\'eralisation de l'identit\'e de Scott sur les
permanents, {\it Linear Algebra Appl.} (to appear)|

\livre \rJou|J.-P. Jouanolou|Polyn\^omes cyclotomiques~: 
Th\'eorie \'el\'ementaire et applica\-tions|Pr\'epublication, Universit\'e de Strasbourg I, 96
pages,  {\oldstyle 1990}|

\article \rKit|R. Kittappa|Proof of a conjecture of 1881 on
permanents|Linear and Multilinear Algebra|10|1981|75--82|

\divers \rKratBI|C.    Krattenthaler|An alternative
evaluation of the Andrews--Burge determinant, in: Mathematical
Essays in Honor of Gian-Carlo Rota, B.~E.~Sagan, R.~P.~Stanley, eds.,
Progress in Math., vol.~161, Birkh\"auser, Boston,
{\oldstyle1999}, pp.~263--270| 

\divers \rKratBN|C.    Krattenthaler|Advanced
determinant calculus, {\sl S\'eminaire Lotharingien Combin.} {\bf42},
``The Andrews Festschrift'' ({\oldstyle1999}), paper B42q, 67~pp| 

\divers \rLasSquare|A. Lascoux|Square-ice Enumeration, 
{\sl S\'eminaire Lotharingien Combin.} {\bf42}, ``The
Andrews Festschrift'' ({\oldstyle 1999}), paper B42p, 15~pp|


\livre \rMac|I. G. Macdonald|Symmetric functions and Hall 
polynomials, {\rm second edition}|Clarendon Press, Oxford, {\oldstyle 1995}|

\livre \rMincPer|H. Minc|Permanents|Encyclopedia of mathematics 
and its applications, vol. 6, Addison-Wesley, Mass., {\oldstyle 1978}|

\article \rMincScott|H. Minc|On a conjecture of R. F. Scott|Linear 
Algebra Appl.|28|1979|141--153|

\livre \rMuirAB|T.    Muir|The theory of
determinants in the historical order of development, {\rm
4~vols}|Macmillan, London, {\oldstyle  1906}--{\oldstyle 1923}|

\article \rSco|R. F. Scott|Mathematical notes|Messenger 
of Math.|10|1881|142--149|

\article \rSvr|D. Svrtan|Proof of Scott's 
conjecture|Proc. Amer. Math. Soc.|87|1983|203--207|

}

\vskip 1cm
\rightline{\quad
\vtop{\halign{#\hfil\cr
{\eightrm I.R.M.A.} and {\eightrm C.N.R.S.}\cr
Universit\'e Louis Pasteur\cr
7, rue Ren\'e-Descartes\cr
F-67084 Strasbourg, France\cr
{\tt guoniu@math.u-strasbg.fr}\cr
}}
\quad
\vtop{\halign{#\hfil\cr
Institut f\"ur Mathematik\cr
Universit\"at Wien\cr
Strudlhofgasse 4\cr
A-1090 Vienna, Austria\cr
{\tt kratt@ap.univie.ac.at}\cr
}}\quad
}

\bye